\documentclass[a4paper,12pt]{amsart}
\usepackage{ifthen}
\usepackage{graphicx}
\usepackage{mathrsfs}
\usepackage{color}
\nonstopmode \numberwithin{equation}{section}
\makeatletter
\@namedef{subjclassname@2020}{\textup{2020} Mathematics Subject Classification}
\makeatother

\newtheorem{thm}{Theorem}[section]
\newtheorem{cor}[thm]{Corollary}
\newtheorem{lem}[thm]{Lemma}
\newtheorem{prop}[thm]{Proposition}

\theoremstyle{definition}

\newenvironment{pf}[1][]{%
 \vskip 3mm
 \noindent
 \ifthenelse{\equal{#1}{}}%
  {{\slshape Proof. }}%
  {{\slshape #1.} }%
 }%
{\qed\bigskip}

\newcounter{alphabet}

\newcommand{\C}{{\mathbb C}}

\newcommand{\uhp}{{\mathbb H}}

\newcommand{\N}{{\mathbb N}}

\newcommand{\R}{{\mathbb R}}

\newcommand{\T}{{\mathcal T}}

\newcommand{\Z}{{\mathbb Z}}

\newcommand{\sphere}{{\widehat{\mathbb C}}}

\renewcommand{\Im}{{\,\operatorname{Im}\,}}

\renewcommand{\Re}{{\,\operatorname{Re}\,}}

\newcommand{\inv}{^{-1}}

\newcommand{\Gauss}{{\null_2F_1}}

\newcommand{\SL}{{\operatorname{SL}}}

\renewcommand{\arg}{\,{\operatorname{arg}\,}}

\newcommand{\aand}{{\quad\text{and}\quad}}

\newcounter{minutes}
\divide\time by 60
\newcounter{hours}
\multiply\time by 60 \addtocounter{minutes}{-\time}
\begin{document}
\bibliographystyle{amsplain}
\title{Hausdorff moment sequences and hypergeometric functions}

\begin{center}
{\tiny \texttt{FILE:~\jobname .tex,
        printed: \number\year-\number\month-\number\day,
        \thehours.\ifnum\theminutes<10{0}\fi\theminutes}
}
\end{center}
\author[T. Sugawa]{Toshiyuki Sugawa}
\address{Graduate School of Information Sciences \\
Tohoku University\\
Aoba-ku, Sendai 980-8579, Japan}
\email{sugawa@math.is.tohoku.ac.jp}
\author[L.-M.~Wang]{Li-Mei Wang}
\address{School of Statistics,
University of International Business and Economics, No.~10, Huixin
Dongjie, Chaoyang District, Beijing 100029, China}
\email{wangmabel@163.com}
\keywords{Continued fraction expansion, totally monotone, Hausdorff moment, Phragm\`en-Lindel\"of principle}
\subjclass[2020]{Primary 33C05; Secondary 30E05}
\begin{abstract}
P\'olya in 1926 showed that the hypergeometric function $F(z)=\null_2F_1(a,b;c;z)$ has a
totally monotone sequence as its coefficients; that is, $F$ is the generating function
of a Hausdorff moment sequence, when $0\le a\le 1$ and $0\le b\le c.$
In this paper, we give a complete characterization of such hypergeometric functions $F$
in terms of complex parameters $a,b,c.$
To this end, we study the class of
general properties of generating functions of Hausdorff moment sequences
and, in particular, we provide a sufficient condition for the class by making use of
a Phragm\`en-Lindel\"of type theorem.
As an application, we give also a necessary and sufficient condition for a 
shifted hypergeometric function to be universally starlike.
\end{abstract}
\thanks{The second author was supported by a grant of University of International Business and Economics (No. 78210418).
}
\maketitle

\section{Introduction}
A sequence of real numbers $c_0, c_1, c_2, \dots$ is called
totally monotone (or completely monotone) if $\Delta^m c_n\ge 0$
for all $m,n\ge 0,$ where $\Delta^m c_n$ are defined inductively in $m$
by $\Delta^0 c_n=c_n$ and $\Delta^{m+1}c_n=\Delta^m c_n-\Delta^m c_{n+1}$
for $n\ge 0.$
A classical theorem of Hausdorff \cite{Haus21} asserts that a sequence
$c_0, c_1, c_2, \dots$ with $c_0=1$ is totally monotone if and only if
there exists a Borel probability measure $\mu$ on the interval $[0,1]$
such that
$$
c_n=\int_0^1 t^n d\mu(t),\quad n=0,1,2,\dots.
$$
Therefore, totally monotone sequences are also called Hausdorff moment sequences.
The generating function of such a sequence is expressed by
$$
F(z)=\sum_{n=0}^\infty c_nz^n=\sum_{n=0}^\infty \int_0^1 t^nz^n d\mu(t)
=\int_{0}^{1}\frac{d\mu(t)}{1-tz}.
$$
By the above integral representation, we observe that the function $F$ is analytically continued
to the slit domain $\Lambda=\C\setminus[1,+\infty)=\{z\ne 1: |\arg(1-z)|<\pi\}.$
We will denote by $\T$ the set of such functions $F$ represented by Borel probability
measures $\mu$ on $[0,1].$
Wall's theorem \cite[Theorems 27.5 and 69.2]{Wall:anal}
tells us that an analytic function $F(z)$ at the origin
belongs to the class $\T$ if and only if $F$ is expressed in the form (sometimes called a g-fraction)
\begin{equation}\label{eq:Wall}
F(z)=\frac1{1-\dfrac{(1-g_0)g_1z}{1-\dfrac{(1-g_1)g_2z}{1-\ddots}}}
\end{equation}
for some sequence $g_0, g_1, g_2, \dots$ of real numbers with $0\le g_n\le 1.$
Moreover, this continued fraction converges uniformly
on any compact subset of $\Lambda.$
The famous Gauss's continued fraction (see \cite[Theorem 89.1]{Wall:anal})
\begin{equation}\label{eq:Gauss}
\frac{\Gauss(a, b+1; c+1; z)}{\Gauss(a, b; c; z)}
=\frac1{1-\dfrac{\dfrac{a(c-b)}{c(c+1)}z}{1-\dfrac{\dfrac{(b+1)(c-a+1)}{(c+1)(c+2)}z}%
{1-\dfrac{\dfrac{(a+1)(c-b+1)}{(c+2)(c+3)}z}{1-\dfrac{\dfrac{(b+2)(c-a+2)}{(c+3)(c+4)}z}%
{1-\ddots}}}}}
\end{equation}
can be viewed as a special case of Wall's theorem by choosing
$$
g_{2k}=\frac{c-a+k}{c+2k}, \quad
g_{2k+1}=\frac{c-b+k}{c+2k+1}
$$
for $k=0,1,2,\dots.$
Here the parameters must satisfy the conditions $0\le a\le c$ and $-1\le b\le c$
so that $0\le g_n\le 1$ for all $n\ge0.$
When $b=0$ one has $\Gauss(a,0;c;z)\equiv 1,$ and hence
\eqref{eq:Gauss} gives a continued fraction expansion
of the function $\Gauss(a,1; c+1;z)$ for $0\le a\le c$
(or, equivalently, for the function $\Gauss(a,1; c;z)$ for $0\le a\le c-1$).
To the best knowledge of the authors, a continued fraction expansion in the form
\eqref{eq:Wall} is not known for the general hypergeometric function
$\Gauss(a,b;c;z)$ ({\it cf.} \cite[\S 15.3]{CJPVW}) except for the case when $b=1.$
On the other hand, P\'olya \cite{Polya} showed that
$\big\{\frac{(a)_n(b)_n}{(c)_nn!}\big\}$ is a totally monotone sequence
if $0<a\le 1$ and $0<b\le c.$
Thus, in this case, $\Gauss(a,b;c;z)$ belongs to $\T$
(see also Theorem 1.6 and Remark 1.7 in K\"ustner \cite{Kustner:2002}).
We remark that the class $\T$ is utilized to characterize the geometric classes of
universally convex and universally starlike functions in \cite{RSS09}.
Subsequently, it was studied by many authors
(see \cite{BRS15, BRS20, LP16, SWW} and references therein).

We organize the present paper as follows.
In Section 2, we investigate basic properties of the functions in $\T$.
In particular, we give a sufficient condition for an analytic function on $\Lambda$ to belong to the class $\T,$
which is suitable to examine hypergeometric functions.
In Section 3, we summarize definitions of hypergeometric functions
together with their fundamental properties.
Then we give a necessary and sufficient condition for the hypergeometric function $\Gauss(a,b;c;z)$
to belong to $\T$ in terms of its complex parameters $a,b,c$
in Theorems \ref{thm:complex} and \ref{thm:main}.
Section 4 will be devoted to proofs of the main results.
In the final section 5, we give a necessary and sufficient condition for a shifted
hypergeometric function to be universally starlike in terms of its parameters.

\section{Preliminaries}

We begin with the following simple facts about totally monotone sequences.

\begin{lem}\label{lem:finite}
Let $c_0=1,c_1,c_2,\dots$ be a totally monotone sequence with
representing probability measure $\mu.$
Then the following assertions hold.
\begin{enumerate}
\item[(i)]
If $c_{m}=c_{m+1}$ for some $m\ge 1,$ then $\mu=(1-k)\delta_0+k\delta_1$
for some $k\in[0,1]$ and $c_1=c_2=c_3=\cdots=k,$ where $\delta_t$ denotes
the Dirac measure with unit mass at $t.$
\item[(ii)]
If $c_m=0$ for some $m\ge 1,$ then $\mu=\delta_0$ and
$c_1=c_2=c_3=\cdots=0.$
\end{enumerate}
\end{lem}

\begin{pf}
Assume that $c_m=c_{m+1}$ for some $m\ge 1.$
Since $\mu$ is a representing measure, we have
$$
c_m-c_{m+1}
=\int_0^1(t^m-t^{m+1})d\mu(t)
=\int_0^1t^m(1-t)d\mu(t)=0.
$$
Hence, the support of the measure $\mu$ is contained in the set $\{0,1\};$
in other words, $\mu=(1-k)\delta_0+k\delta_1,$ where $k=\mu(\{1\}).$
Then
$$
\sum_{n=0}^\infty c_nz^n
=(1-k)\int_0^1\frac{d\delta_0(t)}{1-tz}+k\int_0^1\frac{d\delta_1(t)}{1-tz}
=1-k+\frac{k}{1-z}=1+\frac{kz}{1-z}.
$$
In particular, $c_1=c_2=\cdots=k$ and thus (i) follows.

Assume next that $c_m=0$ for some $m\ge 1.$
Since $c_m\ge c_{m+1}\ge 0,$ we have $c_{m+1}=0.$
Now the conclusion follows from (i).
\end{pf}

We remark that assertion (ii) above appears already in the footnote
of \cite[p.~189]{Polya}.

We recall that the class $\T$ consists of analytic functions $F$ on the slit domain
$\Lambda=\C\setminus[1,+\infty)$ represented by
\begin{equation}\label{eq:mu}
F(z)=F_\mu(z)=\int_{0}^{1}\frac{d\mu(t)}{1-tz},\quad z\in\Lambda,
\end{equation}
for Borel probability measures $\mu$ on the interval $[0,1].$
We summarize basic properties of the functions in $\T$ here.

\begin{prop}\label{prop:T}
\begin{enumerate}
\item The set $\T$ is convex; in other words,
$(1-s)F_0+sF_1\in\T$ whenever $F_0, F_1\in\T$ and $0\le s\le 1.$
\item The set $\T$ is compact with respect to the topology of local uniform convergence.
That is to say, if $f_n\in\T~ (n=1,2,3,\dots)$ converges to a function $f$
locally uniformly on $\Lambda$ as $n\to\infty,$ then $f\in\T.$
\item
For $F\in\T,$ the function $F(x)$ is strictly increasing in $-\infty<x<1$ unless $F\equiv 1.$
\item
If $F=F_\mu\in\T$ for a Borel probability measure $\mu$ on $[0,1],$ then
$$
\lim_{x\to-\infty}F(x)=\mu(\{0\})
\aand
\lim_{x\to 1^-}(1-x)F(x)=\mu(\{1\}).
$$
In particular, the both limits take values in the interval $[0,1].$
\item
If $F\in\T,$ then the function $\hat F(z)=(1-z)\inv F(\frac{z}{z-1})$ is also in $\T.$
\end{enumerate}
\end{prop}

\begin{pf}
The proof of (1) and (2) is easy (see \cite{SWW}).
We next show (3).
Let $F=F_\mu$ for some $\mu.$
Then, for $x<y<1,$
$$
F(y)-F(x)=(y-x)\,\int_0^1\frac{t}{(1-tx)(1-ty)}d\mu(t)\ge 0.
$$
If equality holds, the function $t$ must vanish almost everywhere with respect to
the measure $\mu$ (say $\mu$-a.e.~for short)
and hence, $\mu((0,1])=0,$ in other words, $\mu$ has
a unit mass at $t=0$ and thus $F_\mu=1.$

Finally, we show (4).
For brevity, we set $\alpha_0=\mu(\{0\})$ and $\alpha_1=\mu(\{1\}).$
We also set $\mu_0=\mu-\alpha_0\delta_0-\alpha_1\delta_1,$ where $\delta_t$
stands for the Dirac measure with unit mass at $t.$
Then we can write
\begin{equation}\label{eq:F0}
F(z)=\alpha_0+\frac{\alpha_1}{1-z}+\int_0^1\frac{d\mu_0(t)}{1-tz}
=\alpha_0+\frac{\alpha_1}{1-z}+F_{\mu_0}(z),
\quad z\in\Lambda.
\end{equation}
Since the function $1/(1-tx)$ monotonically tend to 0 as $x\to -\infty$
$\mu_0$-a.e.~on $[0,1],$ we have
$F_{\mu_0}(x)\to 0$ as $x\to-\infty$ by the monotone convergence theorem
(also known as Beppo-Levi's lemma).
Similarly, we have $(1-x)F_{\mu_0}(x)\to 0$ as $x\to 1^-.$
Therefore, it is now easy to derive the assertion (4).
Finally, we show (5).
It follows from Theorem 1.4 in \cite{RSS09} with the choice $\alpha=1.$
However, for convenience of the reader, we provide a direct proof here.
If $F=F_\mu$ for a Borel probability measure $\mu$ on $[0,1],$ we have
$$
\hat F(z)=(1-z)\inv \int_0^1\frac{d\mu(t)}{1-\frac{tz}{z-1}}
=\int_0^1\frac{d\mu(t)}{1-z+tz}
=\int_0^1\frac{d\nu(u)}{1-uz},
$$
where $\nu$ is defined by $d\nu(u)=d\mu(1-u)$ on $[0,1].$
That is, $\hat F=F_\nu\in\T.$
\end{pf}

As is clear by its proof, $\hat{\hat F}=F$ in Proposition \ref{prop:T} (5) and
\begin{equation}\label{eq:dual}
\lim_{x\to-\infty}F(x)=\lim_{t\to 1^-}(1-t)\hat F(t)
\aand
\lim_{x\to 1^-}(1-x)F(x)=\lim_{t\to-\infty}\hat F(t).
\end{equation}

Liu and Pego \cite{LP16} simplified a result in \cite{RSS09} and showed the following
(see Theorem 1 and Remark 2 in \cite{LP16}).

\begin{lem}\label{lem:LP}
Let $F(z)$ be an analytic function on $\Lambda=\mathbb{C}\setminus[1,+\infty)$.
Then $F\in\T$ if and only if the following four conditions are fulfilled:
\begin{enumerate}
\item[(i)] $F(0)=1;$%\medskip
\item[(ii)] $\Im F(z)\ge 0$ for $\Im z>0;$ %\medskip
\item[(iii)] $F(x)\in \R$ for $x\in(-\infty,1);$
\item[(iv)] $\limsup_{x\to -\infty}F(x)\ge 0.$
\end{enumerate}
\end{lem}

When we try to apply this lemma to hypergeometric functions, it is not easy to
check condition (ii).
Instead, we often have enough information about the boundary values of the function $F(z)$
on the slit $(1,+\infty).$
Therefore, we should see that $\Im F(x+i0^+)\ge 0$ for $1<x<+\infty.$
In order to check (ii), however, we need a Phragm\`en-Lindel\"of type theorem in addition to the
information about boundary values.
First we recall the following form of the theorem (see \cite[Theorem 2.3.5]{Rans}).

\begin{thm}[Phragm\`en-Lindel\"of principle]\label{thm:PLP}
Let $S$ be the strip $\{z\in\C: a<\Im z<b\}$ for some real numbers $a<b.$
Suppose that a subharmonic function $U$ on $S$ satisfies the inequality
$\limsup_{z\to\zeta}U(z)\le 0$ for every finite boundary point $\zeta$ of $S.$
If
$$
U(x+iy)\le Ae^{\alpha|x|} \quad (x+iy\in S)
$$
for some constants $0<A<+\infty$ and $\alpha<\pi/(b-a),$
then $U\le 0$ on $S.$
\end{thm}

By using this principle, we can now show the following result.

\begin{lem}\label{lem:PLP}
Let $V$ be a subharmonic function on the upper half-plane $\uhp=\{w: \Im w>0\}.$
Suppose that $V$ satisfies the following three conditions for some constants $A>0$ and $\alpha<1$ and
$\delta\in(0,1):$
\begin{enumerate}
\item[(i)] $\limsup_{\uhp\ni w\to u}V(w)\le 0$ for all $u\in\R\setminus\{0\}$;%\medskip
\item[(ii)] $V(w)\le A|w|^\alpha$ for $w\in\uhp$ with $|w|>1/\delta;$ %\medskip
\item[(iii)] $V(w)\le A|w|^{-\alpha}$ for $w\in\uhp$ with $|w|<\delta.$
\end{enumerate}
Then $V\le 0$ on $\uhp.$
\end{lem}

\begin{pf}
First we observe that for each $u\in\R\setminus\{0\}$ and $\eta>0,$ we have $V\le\eta$
on $\{w\in\uhp: |w-u|\le \varepsilon\}$ for a small enough $\varepsilon>0$ by condition (i).
By a standard compactness argument, we have $V\le A\delta^{\alpha}$ on
$\Omega=\{w\in\uhp: \delta\le |w|\le 1/\delta, \Im w\le \varepsilon \}$
for some $\varepsilon>0.$
By the maximum principle, we have $V\le A\delta^{\alpha}$ on
$\{w\in\uhp: \delta\le |w|\le 1/\delta\}.$
Therefore, $V(w)\le A\max\{|w|^\alpha, |w|^{-\alpha}\}$ on $\uhp.$
Consider now the transformation $w=e^z$ which maps the strip $S=\{z: 0<\Im z<\pi\}$
conformally onto $\uhp.$
Then we apply Theorem \ref{thm:PLP} to the subharmonic function $U(z)=V(e^z)$
to conclude that $U\le 0$ on $S,$ which in turn implies $V\le 0$ on $\uhp$ as required.
\end{pf}

We are now ready to give a sufficient condition for the membership of the class $\T$
which is convenient in the study of hypergeometric functions.

\begin{thm}\label{thm:suff}
Let $F$ be an analytic function on $\Lambda=\C\setminus[1,+\infty)$ with $F(0)=1.$
Suppose that the restriction of $F$ to the upper half-plane $\uhp$ analytically extends
across the slit $(1,+\infty)$ to a function $F^+.$
If the following five conditions are satisfied, then $F\in\T:$
\begin{enumerate}
\item[(1)] $F(x)\in\R$ for all $x\in(-\infty, 1);$
\item[(2)] $\limsup_{x\to-\infty} F(x)\ge 0;$
\item[(3)] $\Im F^+(x)\ge 0$ for all $x\in(1,+\infty);$
\item[(4)] $F(z)=O(|z|^{\alpha})$ as $z\to\infty$ in $\uhp$ for some constant  $\alpha<1;$
\item[(5)] $F(z)=O(|z-1|^{-\alpha})$ as $z\to1$ in $\uhp$ for some constant  $\alpha<1.$
\end{enumerate}
Moreover, $F$ is expressed as in \eqref{eq:mu} for the
Borel probability measure $\mu$ determined by
$$
d\mu(t)=\alpha_0 d\delta_0(t)+\alpha_1 d\delta_1(t)+\frac{\Im F^+(1/t)}{\pi t}dt,
$$
where $\delta_0$ and $\delta_1$ are the Dirac measures with unit mass at $t=0$ and $t=1,$
respectively, and
$$
\alpha_0=\lim_{x\to-\infty}F(x)
\aand
\alpha_1=\lim_{x\to 1^-}(1-x)F(x).
$$
\end{thm}

\begin{pf}
We will apply Lemma \ref{lem:LP} to the function $F.$
In order to check condition (ii) in the lemma, it suffices to show that $\Im F\ge 0$
on $\uhp.$
We can check it by applying Lemma \ref{lem:PLP} to the harmonic function $V(w)=-\Im F(w-1).$
We now conclude that $F\in\T$ by Lemma \ref{lem:LP}.
Let $\mu$ be the representing measure of $F$ as in \eqref{eq:mu}
and write $\mu_0=\mu-\alpha_0\delta_0-\alpha_1\delta_1$ as in the proof of Proposition \ref{prop:T} (4).
Note that $\mu_0$ has its support in $(0,1).$
We now define a Borel probability measure $\nu$ on $(1,+\infty)$ by the relation
$d\nu(u)=u\cdot d\mu_0(1/u).$
Then, by the change of variables $t=1/u,$ \eqref{eq:F0} may be written as
\begin{equation}\label{eq:F1}
F(z)=\alpha_0+\frac{\alpha_1}{1-z}+\int_{(0,1)}\frac{d\mu_0(t)}{1-tz}
=\alpha_0+\frac{\alpha_1}{1-z}+\int_{(1,+\infty)}\frac{d\nu(u)}{u-z}.
\end{equation}
Now the Stieltjes Inversion Formula (see \cite[\S 65]{Wall:anal}) yields
$$
\frac1\pi\int_s^t \Im F^+(x)dx
=\frac1\pi\lim_{y\to0^+}\int_s^t \Im F(x+iy)dx
=\frac{\nu([s,t])+\nu((s,t))}2
$$
for $1<s<t<+\infty.$
This implies that the measure $\nu$ is absolutely continuous and has the density
$\Im F^+(u)/\pi$ on $(1,+\infty);$ and therefore, \eqref{eq:F1} may be
expressed by
$$
F(z)=\alpha_0+\frac{\alpha_1}{1-z}+\frac1\pi \int_1^{+\infty}\frac{\Im F^+(u)}{u-z}du
=\alpha_0+\frac{\alpha_1}{1-z}+\frac1\pi \int_0^{1}\frac{\Im F^+(1/t)}{1-tz}\frac{dt}t.
$$
The proof is now complete.
\end{pf}

We end this section by giving an elementary lemma, which will be needed later.

\begin{lem}\label{lem:mob}
Let $p,q,r,s$ be complex numbers  with $ps-qr\ne0$
and consider the function $f(x)=\dfrac{px+q}{rx+s}.$
If $f(n)\in\R$ for all non-negative integers $n,$ then
$t(p,q,r,s)\in \R^4$ for some $t\in\C\setminus\{0\}.$
In particular, if one of $p,q,r,s$ is a nonzero real number,
$(p,q,r,s)\in\R^4.$
\end{lem}

\begin{pf}
We recall the fact that a M\"obius transformation maps circles or lines to circles or lines.
%is determined by the images of three points in the Riemann sphere.
By the assumption, we see that the M\"obius transformation $f$ maps
the extended real line $\widehat\R=\R\cup\{\infty\}$ onto itself.
Hence, $f$ maps the upper half-plane onto either itself or the lower half-plane.
In the former case, as is well known, $f$ is represented by a matrix in $\SL(2,\R)$ and the conclusion
follows.
In the latter case, $-f$ maps the upper half-plane onto itself and the conclusion follows.
\end{pf}

\section{Main results}

We first recall the definition of hypergeometric functions and their fundamental properties.
For complex parameters
%\footnote{Hypergeometric functions are usually defined for complex
%parameters. In this paper, however, we deal with only real parameters.}
$a,b,c$ with $-c\not\in \N_0$ the hypergeometric function is defined by
$$
\Gauss(a,b;c;z)=\sum_{n=0}^{\infty}\frac{(a)_n(b)_n}{(c)_nn!}z^n
$$
on the unit disk $|z|<1.$
Here, $\N_0=\N\cup\{0\},~ \N=\{1,2,3,\dots\}$ and
$(a)_0=1, ~(a)_{n}=a(a+1)\cdots(a+n-1)=\Gamma(a+n)/\Gamma(a)$ for $n\in\N.$
It is also well known that $w=\Gauss(a,b;c;z)$ is a solution to
the hypergeometric differential equation
\begin{equation}\label{eq:hde}
z(1-z)w''+[c-(a+b+1)z]w'-ab\, w=0,
\end{equation}
which has three regular singularities at $z=0,1,\infty$
in the Riemann sphere $\sphere=\C\cup\{\infty\}.$
Therefore, $\Gauss(a,b;c;z)$ can be analytically continued along any curve in $\C\setminus\{0,1\}.$
In particular, it is well-defined on the slit domain $\Lambda$ as a single-valued function,
and it analytically extends across the interval $(1,+\infty).$
%By the form, we note that $\Gauss(a,b;c;z)$ is symmetric in the parameters $a$ and $b;$
%namely, $\Gauss(a,b;c;z)=\Gauss(b,a;c;z).$
%Therefore, we will always assume that $a\le b$ unless otherwise stated.
Note that $\Gauss(a,b;c;z)\equiv 1$ if $ab=0.$
We often assume that $ab\ne 0$ to avoid this trivial case.

In this paper, our main concern is on the problem asking when $F(z)=\Gauss(a,b;c;z)$
belongs to the class $\T.$
As a test case, it is illustrative to consider the function
$F(z)=\Gauss(a,b;b;z)=(1-z)^{-a}.$
Assume that $F\in\T.$
Then, by Proposition \ref{prop:T}, $F(-\infty)$ should be finite, which enforces $a\ge 0.$
Also, $(1-x)F(x)=(1-x)^{1-a}$ should have a finite limit as $x\to1^-,$
which enforces $a\le 1.$
Note that these limits are 0 for $0<a<1.$
Now we have the necessary condition $0\le a\le 1.$
Conversely, if we assume $0\le a\le 1,$ then it is easy to check $F\in\T$
by Lemma \ref{lem:LP}. %or by Theorem \ref{thm:suff}.
Hence, we conclude that $(1-z)^{-a}$ belongs to the class $\T$ precisely when
$0\le a\le 1.$

The following identity known as Pfaff's formula is useful (see \cite[15.3.4]{AS}
or \cite[p.~79]{AAR}):
\begin{equation}\label{eq:id}
\Gauss(a,b;c;z)=(1-z)^{-a}\Gauss\left(a,c-b;c;\tfrac{z}{z-1}\right).
\end{equation}
In what follows, we will say that the parameters $a,b,c$ are \emph{generic} if none of
$c, a-b, a+b-c$ is an integer\footnote{
More precisely, the parameters are generic if differences of two roots of the indicial equation
of the differential equation \eqref{eq:hde} are not integers at the regular singular points
$z=0,1,\infty.$
}.
%To avoid trivial case, we will always assume that $ab\ne 0.$
We need the following transformation formula, which is valid for generic parameters
(see \cite[15.3.7]{AS}):
\begin{equation}\label{eq:tr}
\begin{aligned}
\Gauss(a,b;c;z)
&=
\frac{\Gamma(c)\Gamma(b-a)}{\Gamma(b)\Gamma(c-a)}(-z)^{-a}
\Gauss\left(a,1-c+a;1-b+a;\frac{1}{z}\right) \\
&+\frac{\Gamma(c)\Gamma(a-b)}{\Gamma(a)\Gamma(c-b)}
(-z)^{-b}\Gauss\left(b,1-c+b;1-a+b;\frac{1}{z}\right)
\end{aligned}
\end{equation}
on $|\arg(-z)|<\pi,$ where the branch of powers of $-z$ should be taken
so that they assume positive values when $z$ is a negative real number.
Note here that \eqref{eq:tr} is valid also for the non-generic case when $c-a-b\in\Z$
as long as $b-a \not\in\Z.$
When $b-a\in\N_0,$ we need the following formula (see \cite[15.3.13, 15.3.14]{AS}):
%\begin{equation}\label{eq:aa}
%\Gauss(a,a;c;z)=\frac{\Gamma(c)}{\Gamma(a)\Gamma(c-a)}(-z)^{-a}
%\sum_{n=0}^{\infty}\frac{(a)_n(1-c+a)_n}{(n!)^2}\big(\log(-z)+C_n\big)(-z)^{-n}
%\end{equation}
%for $|\arg(-z)|<\pi, |z|>1$ and $c-a\not\in\Z,$ where
%$$
%C_n=2\psi(n+1)-\psi(a+n)-\psi(c-a-n)
%$$
%and $\psi(x)=\Gamma'(x)/\Gamma(x)$ denotes the digamma function.
\begin{equation}\label{eq:aa}
\begin{aligned}
&\Gauss(a,a+m;c;z)
=\frac{\Gamma(c)}{\Gamma(a+m)}(-z)^{-a}
\sum_{n=0}^{m-1}\frac{(m-n-1)!(a)_n}{n!\Gamma(c-a-n)}z^{-n} \\
&\quad +\frac{\Gamma(c)}{\Gamma(a+m)\Gamma(c-a)}(-z)^{-a-m}
\sum_{n=0}^{\infty}\frac{(a)_{n+m}(1-c+a)_{n+m}}{n!(n+m)!}\big(\log(-z)+C_n\big)z^{-n}
\end{aligned}
\end{equation}
for $|\arg(-z)|<\pi, |z|>1$ and $c-a\not\in\Z,$ where
$$
C_n=\psi(1+m+n)+\psi(1+n)-\psi(a+m+n)-\psi(c-a-m-n)
$$
and $\psi(x)=\Gamma'(x)/\Gamma(x)$ denotes the digamma function.
When $m=0,$ the first summand in the above is defined to be 0.

If we apply \eqref{eq:id} to the two terms in the right-hand side of \eqref{eq:tr},
we obtain the following formula in $|\arg(1-z)|<\pi$ for generic parameters
(see also \cite[15.3.8]{AS}):
\begin{equation}\label{eq:tr2}
\begin{aligned}
\Gauss(a,b;c;z)&=
(1-z)^{-a} \frac{\Gamma(c)\Gamma(b-a)}{\Gamma(b) \Gamma(c-a)}
\Gauss\left(a, c-b; 1+a-b; \frac{1}{1-z}\right)
\\ & \quad+(1-z)^{-b} \frac{\Gamma(c) \Gamma(a-b)}{\Gamma(a) \Gamma(c-b)}
\Gauss\left(b, c-a; 1+b-a; \frac{1}{1-z}\right).
\end{aligned}
\end{equation}

The following asymptotic behaviour of hypergeometric functions as $z\to1$
in the slit domain $\Lambda$ is also useful (see \cite{AS} or \cite{AAR}).

\begin{lem}\label{lem:asymptotic}
Let $a,~b,~c\in\R$ with $c\ne 0,-1,-2,\dots.$
Put $\Delta=\Delta(a,b,c)=c-a-b.$
\begin{enumerate}
\item
If $\Delta>0,$ the limit of $\Gauss(a,b;c;z)$ exists as $z\to 1$ in $\Lambda$
and it is given by
\begin{equation*}%\label{negative}
\Gauss(a,b;c;1^-)=\frac{\Gamma(c)\Gamma(c-a-b)}{\Gamma(c-a)\Gamma(c-b)}.
\end{equation*}
\item
If $\Delta=0,$ as $z\to 1$ in $\Lambda$
\begin{equation*}\label{zero}
\Gauss(a,b;c;z)=\frac{\Gamma(c)}{\Gamma(a)\Gamma(b)}\log\frac1{1-z}+O(1).
\end{equation*}
\item
The case when $\Delta<0$ reduces to $(1)$ with $\Delta(c-a,c-b,c)=a+b-c=-\Delta>0$
by the transformation formula
$$
\Gauss(a,b;c;z)=(1-z)^\Delta \Gauss(c-a,c-b; c; z).
$$
In particular, as $z\to 1$ in $\Lambda$
$$
\Gauss(a,b;c;z)
=\frac{\Gamma(c)\Gamma(a+b-c)}{\Gamma(a)\Gamma(b)}(1-z)^\Delta+o(|1-z|^\Delta).
$$
\end{enumerate}
\end{lem}

%We are now in a position to state our main results.
Our main objective in this paper is to find a necessary and sufficient condition
for $\Gauss(a,b;c;z)$ to be a member of $\T$.
One may expect that $\Gauss(a,\bar a;c;z)$ belongs to $\T$ for some non-real $a$
and $c>0.$
As an example, consider the case when $a=\bar b=(2+i)/4, c=2.$
The authors confirmed by using Mathematica that
the sequence $\gamma_n=|(\frac{2+i}4)_n|^2/[(2)_n(1)_n]$ satisfies
$\Delta^m\gamma_n\ge 0$ for $m+n\le 1000.$
However, rather surprisingly,  $\Gauss(a,\bar a;c;z)$ cannot belong to $\T$
except for the trivial case:
$(1-z)^{-k}$ for some $k\in[0,1];$ in this case, $(a,b,c)=(k,c,c)$ or $(c,k,c).$
We first give a necessary condition for the membership of $\T.$

\begin{thm}\label{thm:complex}
Let $a, b$ and $c$ be complex numbers with $c\ne 0,-1,-2,\dots.$
Suppose that the hypergeometric function $F(z)=\Gauss(a,b;c;z)$ belongs to $\T.$
Then either $a,b,c\in\R$ or else $F(z)$ is of the form $(1-z)^{-k}$ for a constant
$0\le k\le 1.$
\end{thm}

By virtue of Theorem \ref{thm:complex},
we only need to consider real parameters $a,b,c$ for our main aim.
Note that $\Gauss(a,b;c;z)$ is symmetric in the parameters $a$ and $b;$
namely, $\Gauss(a,b;c;z)=\Gauss(b,a;c;z).$
Therefore, we will always assume that $a\le b$ in the following unless otherwise stated.
We recall that P\'olya \cite{Polya} essentially showed that $\Gauss(a,b;c;z)$ is in $\T$
if $0<a\le 1$ and $0<b\le c$ (without the assumption $a\le b$).
We have the following characterization of $\Gauss(a,b;c;z)$ contained in $\T.$

\begin{thm}\label{thm:main}
Let $a,b,c$ be real numbers with $a\le b,~ ab\ne 0$ and $c\ne 0,-1,-2,\dots.$
The hypergeometric function $\Gauss(a,b;c;z)$ is in the class $\T$
if and only if $0<a\le 1$ and $c\geq a+\max\{0,b-1\}.$
\end{thm}

This result can be rephrased in terms of the coefficients of hypergeometric series.
It complements P\'olya's original result in \cite{Polya}.

\begin{cor}
Let $a,b,c$ be real numbers with $a\le b,~ ab\ne 0$ and $c\ne 0,-1,-2,\dots.$
Then the sequence
$$
\frac{(a)_n(b)_n}{(c)_n(1)_n},\quad n=0,1,2,\dots,
$$
is totally monotone if and only if $0<a\le 1$ and $c\geq a+\max\{0,b-1\}.$
\end{cor}

As for the representing measure of $F(z)=\Gauss(a,b;c;z)$ when $F\in\T,$
we have the following expressions.

\begin{thm}\label{thm:main2}
Assume that real numbers $a,b,c$ satisfy the conditions
$0<a\le \min\{b,1\}$ and $a+\max\{0,b-1\}\le c.$
Then
$$
\Gauss(a,b;c;z)=\frac{\alpha_1}{1-z}+\frac1\pi\int_0^1\frac{h(t)}{1-tz}dt,
\quad z\in\Lambda.
$$
Here
$$
\alpha_1=\begin{cases}
\dfrac{\Gamma(c)}{\Gamma(a)\Gamma(b)}=\dfrac{\Gamma(a+b-1)}{\Gamma(a)\Gamma(b)}
\smallskip
& ~\text{if}~ a+b=c+1, \\
0 & ~\text{otherwise}
\end{cases}
$$
and, when $b-a$ is not an integer,
$$
h(t)=At^{a-1}\Gauss(a,1-c+a;1-b+a;t)+Bt^{b-1}\Gauss(b,1-c+b;1-a+b;t),
$$
where
$$
A=\frac{\Gamma(c)\Gamma(b-a)}{\Gamma(b)\Gamma(c-a)}\sin a\pi
\aand
B=\frac{\Gamma(c)\Gamma(a-b)}{\Gamma(a)\Gamma(c-b)}\sin b\pi,
$$
and, when $m=b-a\in\N_0,$
\begin{align}\label{eq:non-generic}
h(t)
&=\frac{\Gamma(c)}{\Gamma(a+m)}\sin a\pi
\sum_{n=0}^{m-1}\frac{(m-n-1)!(a)_n}{n!\Gamma(c-a-n)}t^{a+n-1} \\
\notag
&+\frac{(-1)^m\Gamma(c)}{\Gamma(a+m)\Gamma(c-a)}
\sum_{n=0}^{\infty}\frac{(a)_{n+m}(1-c+a)_{n+m}}{n!(n+m)!}
\big((C_n-\log t)\sin a\pi-\pi\cos a\pi\big)t^{a+m+n-1}.
\end{align}
\end{thm}

In the above statement, we define $A=0$ when $c-a$
is a non-positive integer. We understand $B$ likewise.

\section{Proof of the main results}

\begin{pf}[Proof of Theorem \ref{thm:complex}]
Suppose $F(z)=\Gauss(a,b;c;z)\in\T$
for some complex parameters $a, b$ and $c$.
We can assume that $ab\ne0$ to avoid the trivial case.
Put $\gamma_n=\frac{(a)_n(b)_n}{(c)_nn!}$ for brevity.
Then, by assumption, the sequence $\gamma_0,\gamma_1,\gamma_2,\dots$ is
totally monotone.
By Lemma \ref{lem:finite}, we have $\gamma_n>0$ for all $n.$
We also have $0<\gamma_{n+1}\le\gamma_n$ and thus
$$
\frac{\gamma_{n+1}}{\gamma_n}=\frac{(a+n)(b+n)}{(c+n)(1+n)}\in(0,1]
$$
for $n\in\N_0.$
Hence,
$$
\frac{(a+n)(b+n)}{c+n}=n+\frac{(a+b-c)n+ab}{n+c}\in\R
$$
for $n\in\N_0.$
We now compute
$$
\det\begin{pmatrix}a+b-c & ab \\ 1&c\end{pmatrix}
=(a+b-c)c-ab=-(a-c)(b-c).
$$
If $a=c$ or $b=c,$ then $F(z)=(1-z)^{-b}$ or $F(z)=(1-z)^{-a},$ respectively.
Thus $F(z)$ is of the form $(1-z)^{-k}$ for some $0\le k\le 1$ as we saw above.
We now assume that $(a-c)(b-c)\ne 0.$
Then Lemma \ref{lem:mob} implies that $a+b, ab, c\in\R.$
Therefore, either $a$ and $b$ are both real or they are conjugate to each other.
We suppose now that $a=\bar b=\alpha+i\beta\not\in\R.$
Applying \eqref{eq:tr2} for $z=-t<0,$ we obtain
\begin{eqnarray*}
\Gauss(a,\bar{a};c;-t)
&=&(1+t)^{-a}
\frac{\Gamma(c) \Gamma(\bar{a}-a)}{\Gamma(c-a)\Gamma(\bar{a})}
\Gauss\left(a,c-\bar{a};1+a-\bar{a}; \frac{1}{t+1}\right)\\
&\quad & +(1+t)^{-\bar{a}}\frac{\Gamma(c) \Gamma(a-\bar{a})}%
{\Gamma(a)\Gamma(c-\bar{a})}\Gauss\left( \bar{a},c-a;1+\bar{a}-a;\frac{1}{t+1}\right).
\end{eqnarray*}
Thus if $t>0$ is large enough, we find that
\begin{eqnarray*}
\Gauss(a,\bar{a};c;-t)&=&
\frac{2}{(1+t)^{\alpha}}
\Re\left[e^{-i\beta\log(1+t)}\frac{\Gamma(c) \Gamma(\bar{a}-a)}{\Gamma(c-a)\Gamma(\bar{a})}
+O(t\inv)\right]
\end{eqnarray*}
as $t\to+\infty.$
By the last expression, we observe that $\Gauss(a,\bar{a};c;-t)$ changes its sign
infinitely many times as $t\to+\infty.$
Therefore, $\Gauss(a,\bar{a};c;z)$ cannot be a member of $\T$
by Proposition \ref{prop:T} (3).
Now the proof is complete.
\end{pf}

In the case $a=\bar b=(2+i)/4, c=2$ which was observed in the previous section,
Mathematica suggests that $\Gauss(a,\bar{a};c;x)>0$ for $-130,000<x<1$
but it assumes a negative value when $x=-131,000.$

\begin{pf}[Proof of Theorem \ref{thm:main}]
Recall that $a,b,c$ are real numbers with $a\le b,~ ab\ne 0$ and $c\ne 0,-1,-2,\dots.$
Put $F(z)=\Gauss(a,b;c;z)$ for brevity and set
$$
\hat F(z)=(1-z)\inv F\left(\tfrac z{z-1}\right)
=(1-z)^{a-1}\Gauss(a,c-b;c; z),
$$
where we have used the formula \eqref{eq:id}.
By using \eqref{eq:tr}, for generic parameters, we have
\begin{equation}\label{eq:ImF}
\Im F^+(x)=Ax^{-a}\Gauss(a,1-c+a;1-b+a;\tfrac1x)+Bx^{-b}\Gauss(b,1-c+b;1-a+b;\tfrac1x)
\end{equation}
for $1<x<+\infty.$

We first assume that $F\in\T.$
Then $\hat F\in\T$ by Proposition \ref{prop:T} (5).
We now show that the conditions $0<a\le 1$
and $c\geq a+\max\{0,b-1\}$ are necessary.
We compute by using \eqref{eq:dual}
$$
\alpha_0
=\lim_{x\to-\infty}F(x)
=\lim_{t\to 1^-}(1-t)\hat F(t)
=\lim_{t\to 1^-}(1-t)^{a}\Gauss(a,c-b;c;t).
$$
According to Proposition \ref{prop:T} (4), $0\le\alpha_0\le 1.$
In view of Lemma \ref{lem:asymptotic} together with $\Delta=c-a-(c-b)=b-a\ge 0,$
the condition $a>0$ is necessary.
We note that $\alpha_0=0$ in this case.
Since $\hat F(x)$ is increasing in $-\infty<x<1$ by Proposition \ref{prop:T} (3),
we have
$$%\begin{align*}
1=\hat F(0)<\lim_{x\to 1^-}\hat F(x)
=\lim_{x\to 1^-}(1-x)^{a-1}\Gauss(a,c-b;c; x).
$$%\end{align*}
By Lemma \ref{lem:asymptotic}, the condition $a-1\le 0$ must be satisfied.
Thus, we have shown that $0<a\le 1.$
By Lemma \ref{lem:asymptotic} again,
$$
\alpha_1=\lim_{x\to1^-}(1-x)\Gauss(a,b;c;x)
$$
is finite if and only if $\Delta=c-a-b\ge -1.$
Hence the inequality $c\ge a+b-1$ should hold by Proposition \ref{prop:T} (3).
Note that $\alpha_1>0$ only if $c=a+b-1.$
In particular, $\alpha_1=0$ for generic parameters.
Finally, we show the inequality $c\ge a.$
It is enough to show it when $b<1.$
Since $F'(0)=ab/c>0,$ we have $c>0.$
For generic parameters, by the above expression of $\Im F^+(x),$ we have
$$
0\le \lim_{x\to+\infty}x^a \Im F^+(x)=A
=\frac{\Gamma(c)\Gamma(b-a)}{\Gamma(b)\Gamma(c-a)}\sin a\pi.
$$
Since $c>0, b>a>0, 0<a<1$ for the generic case, the inequality
$\Gamma(c-a)>0$ holds so that $c>a.$
Note that the same conclusion holds even when $c-a-b\in\Z$ as is noted
right after \eqref{eq:tr}.
We now consider the remaining non-generic case when $b-a\in\N_0$
(in the possible other cases $c\in\N, c-a\in\N_0$ and $c-b\in\N_0,$
we have nothing to prove).
First we show $a\le c$ when $m=b-a\in\N_0.$
If $c-a\in\Z,$ then $c-a\ge 0$ because $c>0$ and $a\le 1.$
Thus we may assume that $c-a\not\in\Z.$
Then, by \eqref{eq:aa}, we obtain
\begin{align}\label{eq:Im}
\notag
&\Im F^+(x)=\Im\left[
\frac{\Gamma(c)}{\Gamma(a+m)}x^{-a}e^{a\pi i}
\sum_{n=0}^{m-1}\frac{(m-n-1)!(a)_n}{n!\Gamma(c-a-n)}x^{-n} \right] \\
\notag
&~ +\Im\left[\frac{\Gamma(c)}{\Gamma(a+m)\Gamma(c-a)}x^{-a-m}e^{(a+m)\pi i}
\sum_{n=0}^{\infty}\frac{(a)_{n+m}(1-c+a)_{n+m}}{n!(n+m)!}\big(\log x-\pi i+C_n\big)x^{-n}\right] \\
%\notag
&=\frac{\Gamma(c)}{\Gamma(a+m)}x^{-a}\sin a\pi
\sum_{n=0}^{m-1}\frac{(m-n-1)!(a)_n}{n!\Gamma(c-a-n)}x^{-n} \\
\notag
&+\frac{(-1)^m\Gamma(c)}{\Gamma(a+m)\Gamma(c-a)}
\sum_{n=0}^{\infty}\frac{(a)_{n+m}(1-c+a)_{n+m}}{n!(n+m)!}
\big((\log x+C_n)\sin a\pi-\pi\cos a\pi\big)x^{-a-m-n}
\end{align}
for $1<x<+\infty.$
In particular, when $m=0,$ we see that
$$
0\le \lim_{x\to+\infty}\frac{x^a}{\log x}\Im F^+(x)
=\frac{\Gamma(c)}{\Gamma(a)\Gamma(c-a)}\sin a\pi.
$$
Hence, we conclude that $\Gamma(c-a)>0.$
Since $c-a\ge c-1>-1$ in this case, it implies the inequality $c-a>0$ as required.
When $m\ge 1,$ by \eqref{eq:Im}, we have
$$
0\le \lim_{x\to+\infty}x^a \Im F^+(x)
=\frac{\Gamma(c)(m-1)!}{\Gamma(a+m)\Gamma(c-a)}\sin a\pi.
$$
Therefore, similarly, we have $c-a>0$ in this case, too.
We have now completed the necessity part.

Next we show the sufficiency part.
To this end, by virtue of Proposition \ref{prop:T} (2), it is enough
to show it only for a dense set of parameters.
From now on, we assume that $a,b,c$ are generic parameters satisfying
$0<a<\min\{b,1\}$ and $a+\max\{0,b-1\}<c.$
By Euler's integral representation, when $0<b<c,$
$$
F(z)=\Gauss(a,b;c;z)=\frac{\Gamma(c)}{\Gamma(b)\Gamma(c-b)}
\int_{0}^{1}t^{b-1}(1-t)^{c-b-1}(1-tz)^{-a}dt,
\quad z\in \Lambda.
$$
Since $\arg[(1-tz)^{-a}]=-a\arg(1-tz)=a\arg(1+tz)\in(0,\pi)$
for $z\in\uhp$ and $0<t<1,$ we see that $(1-tz)^{-a}$
maps the upper half-plane $\uhp$ into itself.
Therefore, by Lemma \ref{lem:LP}, it is easy to see that $F\in\T$ in this case.
We now assume that $c<b.$
If $b\le 1,$ then we apply Euler's integral representation to $\Gauss(b,a;c;z)$
to see that $F\in\T.$
Therefore, we can further assume $b>1$ so that $a+b-1<c.$
Let
$$
h_1(x)=x^{-a}\Gauss(a,1-c+a;1-b+a;1/x)
\aand
h_2(x)=x^{-b}\Gauss(b,1-c+b;1-a+b;1/x).
$$
Then $\Im F^+(x)=Ah_1(x)+Bh_2(x)$ for $x>1$ by \eqref{eq:ImF}.
Since $b, 1-c+b, 1-a+b$ are all positive, we see that
$h_2(x)>0$ for $x>1.$
Note that $h_1$ and $h_2$ are both solutions to the hypergeometric differential equation
\eqref{eq:hde}.
Therefore, the Wronskian $W=h_1h_2'-h_2h_1'$ has the form
$W(x)=Cx^{-c}(x-1)^{c-a-b-1}$ for a constant $C$ (see \cite[Lemma 3.2.6]{AAR} for example).
By the asymptotic expansion of $W(x)$ as $x\to+\infty,$ we can get $C=a-b.$
Since $(h_1/h_2)'=-W/h_2^2>0$ on $(1,+\infty),$
the inequality
$$
\frac{h_1(x)}{h_2(x)}>\lim_{u\to1^+}\frac{h_1(u)}{h_2(u)}
=\frac{\Gamma(1-b+a)\Gamma(b)\Gamma(1-c+b)}%
{\Gamma(1-a+b)\Gamma(a)\Gamma(1-c+a)}
=:D
$$
holds for each $x\in(1,+\infty),$ where we have made use of Lemma \ref{lem:asymptotic} (3).
By using Euler's reflection formula
\begin{equation*}\label{eq:reflection}
\Gamma(z)\Gamma(1-z)=\frac{\pi}{\sin\pi z},
\end{equation*}
we further compute
$$
\frac AB\cdot D=\frac{\Gamma(b-a)\Gamma(1-b+a)\Gamma(c-b)\Gamma(1-c+b)\sin a\pi}%
{\Gamma(a-b)\Gamma(1-a+b)\Gamma(c-a)\Gamma(1-c+a)\sin b\pi}
=-\frac{\sin (c-a)\pi\cdot\sin a\pi}{\sin(c-b)\pi\cdot\sin b\pi}.
$$
Hence,
$$
\frac AB\cdot D+1
=\frac{\sin (b-a)\pi\cdot\sin (c-a-b)\pi}{\sin(c-b)\pi\cdot\sin b\pi}
=\frac{\sin (a-b)\pi\cdot\sin (a+b-c)\pi}{\sin(c-b)\pi\cdot\sin b\pi}.
$$
Since $A>0,$ we now have the estimate
\begin{align*}
\frac{h(x)}{h_2(x)}&=A\frac{h_1(x)}{h_2(x)}+B>AD+B=B\left(\frac AB\cdot D+1\right) \\
&=\frac{\Gamma(c)\Gamma(a-b)\sin b\pi}{\Gamma(a)\Gamma(c-b)}\cdot
\frac{\sin (a-b)\pi\cdot\sin (a+b-c)\pi}{\sin(c-b)\pi\cdot\sin b\pi} \\
&=\frac{\Gamma(c)\Gamma(1-c+b)\sin (a+b-c)\pi}%
{\Gamma(a)\Gamma(1-a+b)}.
\end{align*}
Here, we again used Euler's reflection formula.
Since $1-c+b>0, a+b-c>0$ and $1-a+b>0,$ we now conclude that $h(x)/h_2(x)>0$
and thus $h(x)=\Im F^+(x)>0$ for $x>1.$
We recall that $F(x)\in\R$ for $-\infty<x<1.$
We now investigate the growth order of $F(z)$ as $z\to\infty$ and $z\to 1$ in $\uhp.$
By \eqref{eq:tr}, we easily see that $F(z)=O(|z|^{-a})=o(1)$ as $z\to \infty.$
In particular, $\lim_{x\to-\infty}F(x)=0.$
Since $\Delta=c-a-b<0$ in this case, we obtain from Lemma \ref{lem:asymptotic} (3)
$F(z)=O(|z-1|^{\Delta})$ as $z\to 1$ in $\Lambda.$
We note that $-\Delta=a+b-c<1$ by assumption.
Thus we are able to apply Theorem \ref{thm:suff} to obtain
$F\in\T$ as required.
\end{pf}

We next show Theorem \ref{thm:main2} though it is almost done in the above proof.

\begin{pf}[Proof of Theorem \ref{thm:main2}]
Let $0<a\le \min\{b,1\}$ and $a+\max\{0,b-1\}\le c.$
We use the symbols in Theorem \ref{thm:suff} for $F(z)=\Gauss(a,b;c;z).$
Recall that $\alpha_0=0$ in this case and $\alpha_1>0$ only when $a+b=c+1.$
When $\Delta=c-a-b=-1,$ by using Lemma \ref{lem:asymptotic} (3), we obtain
$$
\alpha_1=\lim_{x\to1^-}(1-x)\Gauss(a,b;c;x)
=\frac{\Gamma(c)}{\Gamma(a)\Gamma(b)}.
$$
If the parameters are generic, as we saw in \eqref{eq:ImF}, the density function $h(t)$
of the measure $\pi\mu_0$ is given by the formula
\begin{align*}
h(t)&=\frac{\Im F^+(1/t)}{\pi t} \\
&=At^{a-1}\Gauss(a,1-c+a;1-b+a;t)+Bt^{b-1}\Gauss(b,1-c+b;1-a+b;t),
\end{align*}
where
$$
A=\frac{\Gamma(c)\Gamma(b-a)}{\Gamma(b)\Gamma(c-a)}\sin a\pi
\aand
B=\frac{\Gamma(c)\Gamma(a-b)}{\Gamma(a)\Gamma(c-b)}\sin b\pi.
$$
By a limiting process, this formula remains valid for non-generic cases
as long as $A$ and $B$ both stay bounded (possibly to be 0).
Therefore, we need a special care only when $m=b-a\in\N_0.$
In this case, we can also deduce \eqref{eq:non-generic} by using
\eqref{eq:Im}.
\end{pf}

\section{Application to universally starlike functions}

An analytic function $f$ on $\Lambda=\C\setminus[1,+\infty)$ with $f(0)=0, f'(0)=1$ is called
\emph{universally starlike}
if it maps every disk or half-plane $\Omega$ in $\Lambda$ with $0\in\Omega$ univalently
onto a starlike domain with respect to the origin.
An analytic function $f$ on $\Lambda$ with $f(0)=0$ and $f'(0)=1$ is universally starlike
if and only if $zf'(z)/f(z)$ belongs to the class $\T$ (see \cite{RSS09}).
We again examine the function $f(z)=z(1-z)^{-k}$ for a complex constant $k.$
Since $zf'(z)/f(z)=[1+(k-1)z]/(1-z),$ it is easy to check by Lemma \ref{lem:LP}
that $f$ is universally starlike precisely when $k\in[0,1].$
Recall that $f(z)=z\Gauss(k,c;c;z)=z\Gauss(c,k;c;z)$ for an arbitrary $c\in\C\setminus(-\N_0).$
Concerning this case, we also note the following fact.

\begin{lem}\label{lem:us}
Let $f$ be a universally starlike function.
Suppose that $zf'(z)/f(z)$ is represented by a Borel probability measure $\mu$
whose support is contained in the set $\{0,1\}.$
Then $f(z)=z(1-z)^{-k},$ where $k=\mu(\{1\})\in[0,1].$
\end{lem}

\begin{pf}
In view of the proof of Lemma \ref{lem:finite}, we have
$$
\frac{f'(z)}{f(z)}=\frac1z+\frac{k}{1-z}.
$$
We integrate both sides to obtain the relation $\log[f(z)/z]=-k\log(1-z),$
which implies $f(z)=z(1-z)^{-k}$ as required.
\end{pf}

By using K\"ustner's theorem \cite[Theorem 1.5]{Kustner:2002},
it was shown in \cite{RSS09} that the shifted hypergeometric
function $z\Gauss(a,b;c;z)$ is universally starlike if $0<a\le \min\{1,c\}, 0<b\le c.$
It turns out that this condition is also necessary except for the elementary case
$f(z)=z\Gauss(k,c;c;z)=z\Gauss(c,k;c;z).$
Since $\Gauss(a,b;c;z)=\Gauss(b,a;c;z),$ we may assume that $\Re a\le \Re b$ in the sequel
so that $a\le c$ follows from the condition $b\le c.$

\begin{thm}
Let $a,b,c$ be complex parameters with $\Re a\le \Re b$ and $c\ne 0,-1,-2,\dots.$
The shifted hypergeometric function $f(z)=z\Gauss(a,b;c;z)$ is universally starlike
if and only if either $a,b,c$ are positive numbers with $0<a\le 1, 0<b<c$
or $f(z)=z(1-z)^{-k}$ for some $k\in[0,1].$
\end{thm}

\begin{pf}
The sufficiency part follows from K\"ustner's theorem as we mentioned above.
Here we only note the important identity for the function $f(z)=z\Gauss(a,b;c;z)$:
\begin{equation}\label{eq:log-deri}
\frac{zf'(z)}{f(z)}
=1-a+a\frac{\Gauss(a+1,b;c;z)}{\Gauss(a,b;c;z)}.
\end{equation}
K\"ustner \cite{Kustner:2002} indeed showed that the function
$G(z):=\Gauss(a+1,b;c;z)/\Gauss(a,b;c;z)$ belongs to $\T$ under the assumption.

We now show the necessity part.
It is known that $f(z)/z$ belongs to $\T$ whenever $f(z)$ is universally starlike, in general
(see \cite[p.~289]{RSS09}).
Therefore, if $f(z)=z\Gauss(a,b;c;z)$ is universally starlike, then
$F(z):=\Gauss(a,b;c;z)$ belongs to the class $\T.$
We further assume that $f(z)$ is not of the form $z(1-z)^{-k}.$
Then we apply Theorems \ref{thm:complex} and \ref{thm:main} to find
that $a,b,c$ are positive numbers with $0<a\le 1$ and $a+\max\{0,b-1\}\le c.$
Now \eqref{eq:log-deri} implies that $1-a+aG\in\T$.
An elementary computation yields that
\begin{eqnarray*}
1-a+aG
&=&1-a+a\left[1+\frac{b}{c}z+\frac{b((a+b+1)c-ab)}{c^2(c+1)}z^2+\cdots\right]\\
&=&1+\frac{ab}{c}z+\frac{ab((a+b+1)c-ab)}{c^2(c+1)}z^2+\cdots.
\end{eqnarray*}
Since the Taylor coefficients of $1-a+aG(z)$ constitute a totally monotone sequence,
the inequality
$$
\frac{ab}{c}-\frac{ab((a+b+1)c-ab)}{c^2(c+1)}
=\frac{ab[c^2-(a+b)c+ab]}{c^2(c+1)}
=\frac{ab(c-a)(c-b)}{c^2(c+1)}\ge 0
$$
holds.
In view of Lemma \ref{lem:finite} (i) and Lemma \ref{lem:us}, we find that
equality does not hold by assumption.
Hence, we have $(c-a)(c-b)>0$ which implies $b<c$ since $a\le c$.
Now the necessity part has been proved.
\end{pf}

%%%%%%%%%%%%%%%%%%%%%%%%%%%%%%%%%%%%%%%%%%%%%%%%%%%%%%%%%%%%%%%%%%%%%%%%%%%%%%%%%%%%%%%%%%%%%%%%%%%%%%%%%%%%%%%%%%%%%%%%%%%%%%%%%%55

%%%%%%%%%%%%%%%%%%%%%%%%%%%%%%%%%%%%%%%%%%%%%%%%%%%%%%%%%%%%%%%%%%%%%%%%%%%%%%%%%%%%%%%%%%%%%%%%%%%%%%%%%%%%%%%%%%%%%%%%%%%%%%%%%%55
\def\cprime{$'$} \def\cprime{$'$} \def\cprime{$'$}
\providecommand{\bysame}{\leavevmode\hbox
to3em{\hrulefill}\thinspace}
\providecommand{\MR}{\relax\ifhmode\unskip\space\fi MR }
\providecommand{\MRhref}[2]{%
  \href{http://www.ams.org/mathscinet-getitem?mr=#1}{#2}
} \providecommand{\href}[2]{#2}

\end{document}